\documentclass[11pt, letterpaper]{article}

\usepackage[affil-it]{authblk}
\renewcommand{\footnote}{\fnsymbol{footnote}}
\usepackage{epic}
\usepackage{eepic}
\usepackage{graphicx}
\usepackage{algorithm,algorithmic}
\usepackage{tikz}
\usepackage{xcolor,colortbl,ragged2e,rotating}
\usepackage{wrapfig}
\usepackage{amsmath,amsfonts,amssymb,amsthm}
\usepackage{listings}
\usepackage{spverbatim}
\usepackage{subcaption}
\usepackage{hyperref}
\usepackage{mathtools}
\usepackage{lineno}
\usepackage{float}
\usepackage{array}
\usepackage{longtable}
\usepackage{textcomp}
\usepackage{varioref}
\usepackage{xspace}
\usepackage{makeidx}
\usepackage{verbatim}
\usepackage{ifpdf}
\usepackage{tabularx}
\usepackage{keyval}
\usepackage{caption}
\usepackage{cite}
\usepackage{dsfont}
\usepackage{amsmath}
\usepackage{amsfonts}
\usepackage{amssymb}
\usepackage[left=2.00cm, right=2.00cm, top=2.00cm, bottom=2.00cm]{geometry}
\usepackage{varwidth}
\DeclareMathOperator{\sech}{sech}
\graphicspath{{Beam_graphics/}{Beam_tcache/}{Beam_gcache/}}
\DeclareGraphicsExtensions{.pdf,.eps,.ps,.png,.jpg,.jpeg}

\title{\bf{On non-existence of bifurcations in one-dimensional Bratu equation}}

\author[]{Shrinidhi S. Pandurangi\thanks{corresponding author: shrinidhi.pandurangi@iitb.ac.in} }
\author[]{Utkarsh S. Nikam}

\affil[]{Department of Mechanical Engineering, Indian Institute of Technology Bombay, Mumbai, India}

\begin{document}

\maketitle

\begin{abstract}
In this paper, we revisit the classical problem of Bratu differential equation in one-dimension. While it is known that the finite difference discretized form of continuous Bratu equation gives rise to spurious bifurcations, we show that spurious bifurcation points exist even when the finite element approach is employed. We then present an analytical proof demonstrating that there are no bifurcations when the continuous Bratu equation is considered.
\end{abstract}

\section{Introduction}\label{sec:Section_1}

The Bratu (Liouiville-Bratu-Gelfand) equation serves as a model in several applications, including the modeling of thermal ignition \cite{Frank-Kamenetskii} and chemical diffusion processes \cite{Frank-Kamenetskii, JACOBSEN2002283}. It has been used in the theory of chemical reactors \cite{Gavalas1968} and as a model for the temperature distribution of an isothermal gas sphere in gravitational equilibrium \cite{Chandrashekhar}. The Bratu equation is employed as a prototypical example to test the capability of a numerical method to compute multiple solutions to nonlinear parameterized boundary value problems, see \cite{Farrell2015} for example.\\
\\
In this work, we consider the one-dimensional Bratu equation with parameter $\lambda \in \mathds{R}$, 
\begin{equation}
    u''(x) + \lambda e^{u(x)} = 0, \ x\in[-1,1],
    \label{eq.bratu}
\end{equation}
satisfying Dirichlet boundary conditions $u(x=\pm1) = 0$. It is well known that the numerically computed bifurcation diagram for (\ref{eq.bratu}) using the finite difference (FD) method gives rise to spurious bifurcations \cite{govaerts2000numerical, korman2003accurate}. Spurious solutions for the difference equations arising from the discretization are not limited to second order differential equations (ODEs) of the form (\ref{eq.bratu}). Such solutions were obtained by Domokos and Holmes \cite{domokos2003nonlinear} for the ODE that governs beam buckling. Allgower \cite{ALLGOWER19751} has shown that $u'' + \lambda u^k = 0$, $u(0) = u(\pi) = 0$ may admit what they refer to as `numerically irrelevant' solutions (NIS) when the FD approximation is used. More generally, Peitgen \cite{peitgen1980nonlinear} and Beyn and Lorenz \cite{beyn1982spurious} reported the existence of NIS in FD solutions for various boundary value problems of the form $u'' + \lambda f(u) = 0$, where $f$ is an asymptotically linear or superlinear continuous function having several zeros. In particular, they showed that a plethora of spurious solutions are obtained when $f(u) = e^u$. Kuiper \cite{kuiper1988spurious} investigated the discretized version of the FD of (\ref{eq.bratu}) and proved that spurious solutions bifurcating from the branch of symmetric solutions are antisymmetric about $x=0$. Moreover, it was shown that the bifurcations occur only when an odd number of elements are used in discretization. \\ \\
The above discussion regarding the occurrence of spurious bifurcations in the discretized version of the Bratu equation serves as a primary motivation to study the bifurcations of (\ref{eq.bratu}). The objectives of this work are two fold: (1) We report on spurious bifurcations in finite element (FE) approximation to the Bratu equation, which has not been previously reported in the literature to the best of our knowledge.  (2) Using the tools of analytical bifurcation theory, we present a direct proof showing that there are no bifurcations on the non-trivial branch of solutions to the continuous one-dimensional Bratu equation. We emphasize that by bifurcation we necessarily mean branching of solutions\cite{rabinowitz1971, stakgold1971}. Therefore, our result is significant given the rarity of problems in which bifurcations along a non-trivial global solution branch can be ruled out. \\

The paper is organized as follows. In Section \ref{sec:Section_2}, we provide bifurcation diagrams and the locations of bifurcation points using FE and FD-based discretizations. We present our analytical results in Section \ref{sec:Section_3} followed by concluding remarks in Section \ref{sec:Section_4}.

\section{Numerical results}\label{sec:Section_2}

The objective of this section is to show the existence of spurious bifurcation points on the numerically computed bifurcation diagram for (\ref{eq.bratu}). We present these results using the FE method in addition to the FD method reported in the literature (see \cite{govaerts2000numerical}). For each of the two methods, we use a discretization of the domain $x\in[-1, 1]$ with $N$ number of elements such that the grid is uniformly spaced with grid size $h=2/N$ and grid points placed symmetrically about $x=0$. We then employ Newton's method coupled with the numerical continuation technique \cite{keller1987lectures} to solve the resulting parametrized nonlinear algebraic equations. \\

\begin{figure}[H]
    \centering
    \includegraphics[scale = 0.53]{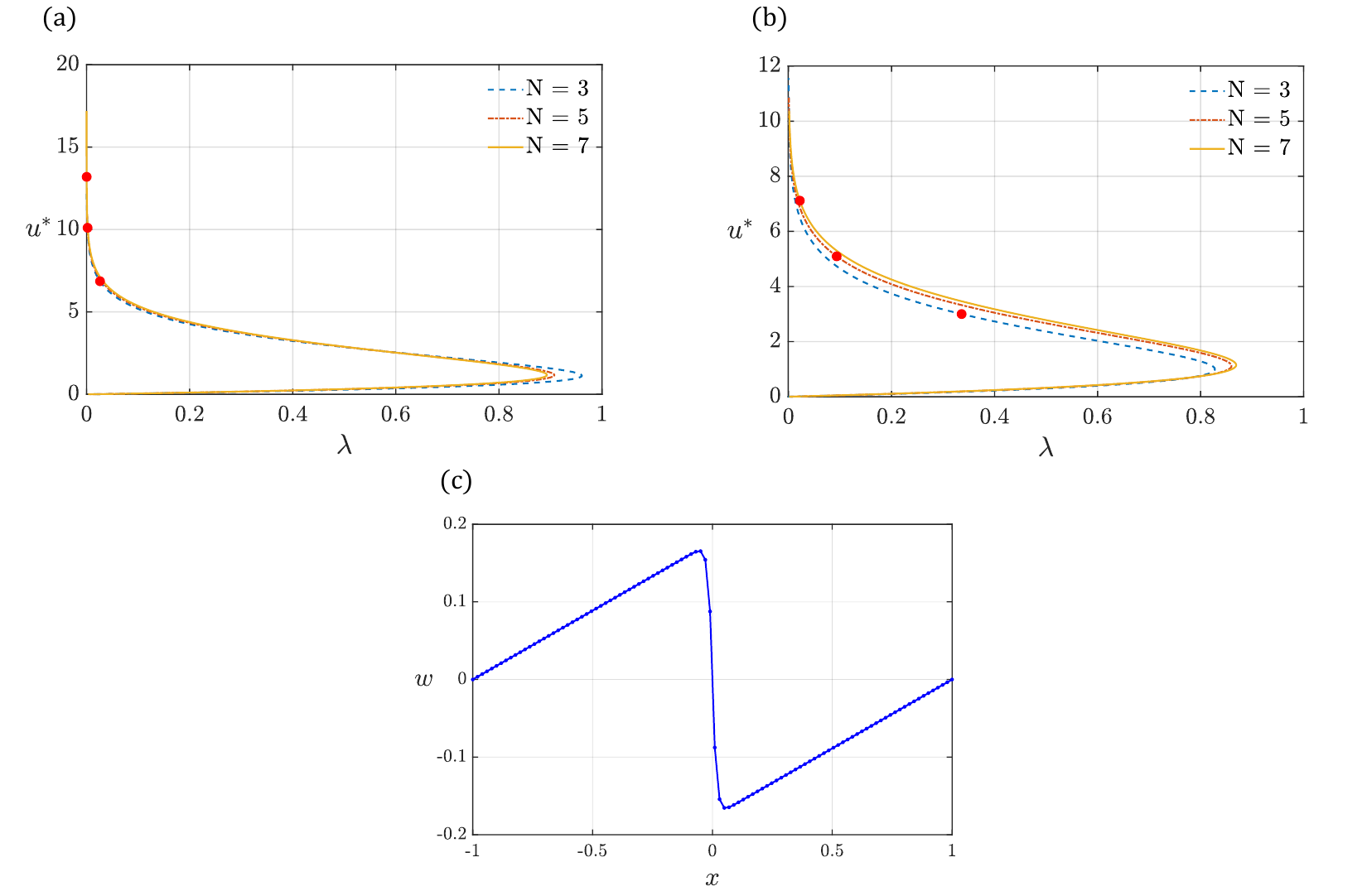}
    \caption{Bifurcation diagram obtained using finite element (a) and finite difference (b) discretizations. The bifurcation points are marked with a red dot on the solution path. (c) Characteristic sawtooth form of the eigenvector $w(x)$ at the bifurcation point for $N = 101$.}
    \label{fig:bif_diag}
\end{figure}

Clearly, $(u(x), \lambda) \equiv (0,0)$ corresponds to the trivial solution of (\ref{eq.bratu}). The non-trivial solutions of (\ref{eq.bratu}), denoted by $(u_0(x), \lambda_0)$, computed using the FE and FD methods are plotted in bifurcation diagrams shown in Fig. \ref{fig:bif_diag}(a) and Fig. \ref{fig:bif_diag}(b) respectively, with the response $u^*:=u(x=0)$ on the vertical axis and the parameter $\lambda$ on the horizontal axis. At every solution point along the bifurcation path, we monitor the eigenvalues of the tangent stiffness matrix. The occurrence of a zero eigenvalue on the path indicates a critical point (a limit point or a bifurcation point). Starting with $(u_0, \lambda_0) = (0,0)$, a zero eigenvalue is first observed at the limit point. Upon continuing along the upper branch, we encounter a second critical point corresponding to a spurious bifurcation point, shown with a red circle for both FE and FD methods. These bifurcation points are observed only when an odd number ($N$) of elements are used in the discretization. The results corresponding to $N = 3, 5, \text{and} \ 7$ are shown in Figures \ref{fig:bif_diag}(a) and \ref{fig:bif_diag}(b). Observe that the bifurcation point moves upwards along increasing $u_0^*$ as the value of $N$ increases.  The exact location of bifurcation points for a few selected number of elements up to $N = 101$ using FE and FD are given in Tables \ref{table1} and \ref{table2} respectively. We conclude that, for both the discretization schemes, the bifurcation point does not converge with increasing numerical refinement, it rather differs significantly for larger values of $N$. This is one of the features of spurious bifurcations, highlighted using FD results in \cite{peitgen1980nonlinear}.

\begin{table}[H]
    \centering
    \begin{tabular}{|c|c|c|c|c|c|}
        \hline
        $N$ & 3 & 5 & 7 & 51 & 101 \\
        \hline
        $\lambda_0$ & 0.0261 & 0.0019 & 1.393$\times$10$^{-4}$ & 7.12$\times$10$^{-27}$ & 1.583$\times$10$^{-49}$ \\
        \hline
        $u_0^*$ & 6.8607 & 10.1049 & 13.1939 & 67.7518 & 121.0915 \\
        \hline
    \end{tabular}
    \vspace{0.2cm}
    \caption{Bifurcation points calculated using the finite-element method.}
    \label{table1}
\end{table}

\begin{table}[h]
    \centering
    \begin{tabular}{|c|c|c|c|c|c|}
        \hline
        $N$ & 3 & 5 & 7 & 51 & 101 \\
        \hline
        $\lambda_0$ & 0.336 & 0.0934 & 0.0219 & 2.42$\times$10$^{-17}$ & 1.147$\times$10$^{-33}$ \\
        \hline
        $u_0^*$ & 3.0002 & 5.0958 & 7.1154 & 45.0789 & 83.9083 \\
        \hline
    \end{tabular}
    \vspace{0.2cm}
    \caption{Bifurcation points calculated using finite-difference method.}
    \label{table2}
\end{table}

Another signature of spurious bifurcations is the lack of reflection symmetry of the eigenvector corresponding to the zero eigenvalue \cite{govaerts2000numerical}. This can be seen in Figure \ref{fig:bif_diag}(c) which shows the eigenvector $w(x)$, having a characteristic ``saw-tooth" pattern antisymmetric about $x=0$, for $N=101$ elements using FE. Similar antisymmetric eigenvectors were obtained for other $N$ values considered in this work, for both FE and FD implementation schemes. Note however that, the numerical solutions $u(x)$ along the non-trivial solution path are symmetric about $x=0$ for all parameter values. 
\\
\textbf{Remark.} While the exact location of spurious bifurcation points differs (see Tables \ref{table1} and \ref{table2}), the occurrence of spurious numerical bifurcation points is not inherent to the choice of the numerical scheme.

\section{Analysis of critical points}\label{sec:Section_3}

In this section, we present an analytical proof demonstrating the non-existence of bifurcation points by directly working with the continuous Bratu equation. In particular, we suggest a change of variable that renders the variable coefficient ODE, resulting from the linearization of (\ref{eq.bratu}) about the non-trivial solution, amenable for analysis. To this end, we begin by recalling the well-known non-trivial solutions of (\ref{eq.bratu}) \cite{liouville1853equation}. They can be written in terms of $u_0^*\in(0,\infty)$ as,
\begin{equation}
\begin{aligned}
        &u_0(x) = -2\ln\left[\exp\left(-u_0^{*}/2\right)\cosh(x\cosh^{-1}\left(\exp\left(u_0^{*}/2\right)\right))\right],
        \\  &\lambda_0 = 2 \left[ \exp\left(-u_0^{*}/2\right)\cosh^{-1}\left(\exp({u_0^{*}/2})\right) \right]^2.
        \label{eq.u_0_star}
\end{aligned}
\end{equation}
\\
To proceed further, we introduce a parameter $\alpha$ such that,
\begin{equation}
    \alpha:=\cosh^{-1}\left(\exp({u_0^{*}/2})\right),
\end{equation}
and rewrite (\ref{eq.u_0_star}) in parameterized form as,
\begin{equation}
    u_0(x) = -2\ln\left[\frac{\cosh(x\alpha)}{\cosh(\alpha)}\right] \ \text{and} \ \lambda_0 = 2 \left[ \frac{\alpha}{\cosh(\alpha)} \right]^2.
    \label{eq.u_0}
\end{equation}
\\
The curve generated through (\ref{eq.u_0}) for $\alpha \in [0, \infty)$ will be referred to as the primary path of solutions to (\ref{eq.bratu}), with $\alpha=0$ corresponding to the trivial solution.\\
\\
Equation (\ref{eq.bratu}) can be written in operator form,
\begin{equation}
    F(u,\lambda) := u'' + \lambda e^u = 0
    \label{eq.bratu_operator}
\end{equation}
with boundary conditions
\begin{equation}
    u(x=-1)=u(x=1)=0,
    \label{eq.bratu_bcs}
\end{equation}
yielding $F(u_0(x), \lambda_0) \equiv 0$. We linearize (\ref{eq.bratu_operator}) about $(u_0(x),\lambda_0)$ and define

\begin{equation}
  \mathcal{L}w(x):=D_{u}F(u_0(x),\lambda_0)w = \left[\frac{d}{d\varepsilon} F(u_0 + \varepsilon w, \lambda_0)\right]_{\varepsilon = 0},
\end{equation}
such that,
\begin{equation}
  \mathcal{L}w(x)= w'' +\lambda_0 e^{u_0(x)}w=0, \label{eq.lin_eq}
\end{equation} 
subject to the boundary conditions,
\begin{equation}
    w(x=-1) = w(x=1) = 0.
    \label{eq.lin_bcs}
\end{equation}
Equation (\ref{eq.lin_eq}) can be rewritten in terms of $\alpha$ as,
\begin{equation}
    \mathcal{L}w(x) = w'' + 2\left[\frac{\alpha}{\cosh(\alpha x)}\right]^2 w = 0.
    \label{eq.lin_eq_alpha}
\end{equation}
When $\alpha=0$, we get $w(x)\equiv0$ as the only solution to (\ref{eq.lin_eq_alpha}) implying that $(u(x), \lambda) \equiv (0,0)$ is a regular solution point. In the following discussion, we will consider the case where $\alpha>0$. We shall refer to a solution point on the primary path as a critical point if there exist non-trivial solutions to (\ref{eq.lin_eq_alpha}) for $\alpha>0$.
\\ \\
\textbf{Lemma 3.1.} For $\alpha>0$, the kernel $\mathcal{N}(\mathcal{L}w(x))$ of the linear operator $\mathcal{L}$ defined in (\ref{eq.lin_eq_alpha}) is \\
(i) nontrivial when $\alpha$ satisfies $\alpha \tanh(\alpha) = 1$,\\
(ii) spanned by $w(x)=\alpha x \tanh(\alpha x) - 1$,\\
(iii) and has $\text{dim} \ \mathcal{N}(\mathcal{L}w(x)))=1$.\\ \\
\textbf{Proof.}
We introduce a change of a variable for $x$ such that $t = \tanh(\alpha x)$.
Upon writing (\ref{eq.lin_eq_alpha}) in terms of the independent variable $t$, we get
\begin{equation}
    (1-t^2)\frac{d^2w(t)}{dt^2} - 2t\frac{dw(t)}{dt} +2w(t) = 0,
    \label{eq.legendrediff}
\end{equation}
subject to the boundary conditions $w(\tanh(\alpha)) = w(-\tanh(\alpha)) = 0$. Note that the parameter $\alpha$ now appears in the boundary conditions.
\\ \\
The equation (\ref{eq.legendrediff}) is a Legendre differential equation of order 1. The general solution of (\ref{eq.legendrediff}) is given by,
\begin{equation}
    w(t) = c_{1}t + c_{2}\left[\frac{-t}{2}\ln\left(\frac{1-t}{1+t}\right)-1\right].
    \label{eq.legendre_gen_sol_t}
\end{equation}
Upon substituting $t$ with $x$ and simplifying (\ref{eq.legendre_gen_sol_t}) we get
\begin{equation}
    w(x) = c_{1}\tanh(\alpha x) + c_2[\alpha x \tanh(\alpha x) - 1].
    \label{eq.legendre_gen_sol}
\end{equation}
\\
The general solution (\ref{eq.legendre_gen_sol}) should satisfy the boundary conditions (\ref{eq.lin_bcs}) which yield two equations for the unknown constants $c_{1}$ and $c_{2}$.
\begin{equation}
\begin{aligned}
    -c_{1}\tanh(\alpha) + c_{2}[\alpha \tanh(\alpha) -1] &= 0,\\
    c_{1}\tanh(\alpha) + c_{2}[\alpha \tanh(\alpha) -1] &= 0. \label{eq.const_eq}
\end{aligned}
\end{equation}
\\
Satisfaction of (\ref{eq.const_eq}) requires that $c_{1}\tanh{\alpha}=0$. As $\alpha > 0$, we get $c_1 = 0$. For nontrivial solutions, the boundary conditions on $w(x)$ require that,
\begin{equation}
    \alpha \tanh(\alpha) = 1.
    \label{eq.criticality_condition}
\end{equation}
Therefore, the kernel of (\ref{eq.lin_eq_alpha}) is spanned by,
\begin{equation}
    w(x) = \alpha x \tanh(\alpha x) - 1,
    \label{eq.particular_sol}
\end{equation}
such that $\alpha>0$ satisfies (\ref{eq.criticality_condition}). \\ \\
Define $f(\alpha):=\alpha \tanh(\alpha)-1$. To prove (iii), it is sufficient to show that the multiplicity of $f^{-1}(0)$ is one. The function $f(\alpha)$ is continuous at every $\alpha \in \mathds{R}$. As $f(\alpha \searrow 0) = -1$ and $f(\alpha \rightarrow \infty) \rightarrow \infty$, by the intermediate value theorem, there exists at least one $\alpha > 0$ such that $f(\alpha) = 0$. Furthermore, as $f(\alpha)$ increases strictly on $(0, \infty)$, it can be concluded that $f(\alpha)= \alpha \tanh (\alpha) - 1$ has exactly one solution for $\alpha>0$ with a multiplicity of one.\qed 
\\ \\
We denote the unique solution of (\ref{eq.criticality_condition}) with $\bar{\alpha}$ and the corresponding $(u(x), \lambda)$ using (\ref{eq.u_0}) as $\bar{u}_0:=u_0(x,\bar{\alpha}), \bar{\lambda}_{0}:=\lambda_{0}(\bar{\alpha})$. \\
\textbf{Remark.} Note that, maximizing $\lambda_{0}$ given in (\ref{eq.u_0}) with respect to $\alpha$ gives the condition (\ref{eq.criticality_condition}) for the limit (turning) point alone. However, we have shown that (\ref{eq.criticality_condition}) is indeed a necessary condition for the existence of all critical points, i.e. bifurcation points and limit points in general. We now proceed to show that the critical point obtained using \ref{eq.criticality_condition} is a limit point, ruling out any bifurcations (branching) on the primary path of solutions \ref{eq.u_0}.\\
\\
It is straightforward to show that the linear operator (\ref{eq.lin_eq_alpha}) is self-adjoint. Therefore, if $\mathcal{L}^{\dagger}$ is the adjoint operator of $\mathcal{L}$, then $w(x)$ given in (\ref{eq.particular_sol}) is a solution to $\mathcal{L}^{\dagger}w(x)=0$. \\
\\
\textbf{Theorem 3.1.} Let (\ref{eq.bratu_operator}) define a continuously differentiable map from all real numbers $\lambda$ and twice-continuously differentiable functions $u(x)$ on $[-1,1]$ satisfying (\ref{eq.bratu_bcs}). If $F(\bar{u}_0, \bar{\lambda}_0)=0$ for some $(u_0,\lambda_0)$ and $\text{dim} \ \mathcal{N}(\mathcal{L}w(x)))=1$ and 
\begin{equation}
    D_{\lambda}F(\bar{u}_0, \bar{\lambda}_0) \notin \mathcal{R}(D_{u}F(\bar{u}_0, \bar{\lambda}_0)),
    \label{eq.limit_pt_cond}
\end{equation}
then there exists a unique continuously differentiable curve $((u_0(s),\lambda_0(s))|s\in(-\delta, \delta)), (u_0(0),\lambda_0(0)) = (\bar{u}_0,\bar{\lambda}_0)$ such that $F(u_0(s), \lambda_0(s))\equiv 0$ \cite{kielhofer2011}.
\\ \\
\textbf{Proof.} As $(\bar{u}_0, \bar{\lambda}_0)$ belongs to the primary path solutions (\ref{eq.u_0}), $F(\bar{u}_0, \bar{\lambda}_0)= 0$ and we have $\text{dim} \ \mathcal{N}(\mathcal{L}w(x)))=1$ from Lemma 3.1. To show (\ref{eq.limit_pt_cond}), we first define the inner product $\langle g(x), h(x) \rangle := \int_{-1}^{1} g(x) h(x) \ dx$. Then (\ref{eq.limit_pt_cond}) is equivalent to showing
\begin{equation}
    \langle w(x), D_{\lambda}F(\bar{u}_0, \bar{\lambda}_0) \rangle \neq 0.
    \label{eq.limit_pt_equiv_cond}
\end{equation}
\\
Upon evaluating the inner product in (\ref{eq.limit_pt_equiv_cond}), we get
\begin{equation}
    \langle w(x), e^{\bar{u}_0(x)}\rangle = \frac{1}{\sech^{2}(\bar{\alpha})}\int_{-1}^{1}\frac{\bar{\alpha} x \tanh(\bar{\alpha} x)-1}{\cosh^{2}(\bar{\alpha} x)}dx=-\left(1+\frac{\sinh(\bar{\alpha})\cosh(\bar{\alpha})}{\bar{\alpha}}\right) \neq 0,
\end{equation}
\\
as $\bar{\alpha}>0$. From the implicit function theorem, we conclude that there exists a continuously differentiable curve $((u_0(s), \lambda_0(s)), s\in(-\delta, \delta))$ corresponding to all solutions of $F(u, \lambda)= 0$ in the neighbourhood of $(u_0(0), \lambda_0(0))=(\bar{u}_0, \bar{\lambda}_0)$ . \qed \\
Theorem 3.1 shows that the only critical point on the primary path is a limit (turning) point.

\section{Concluding remarks}\label{sec:Section_4}

In conclusion, we highlight that despite the choice of discretization that respects the reflection symmetry required by (\ref{eq.bratu}), we get antisymmetric eigenmodes at the spurious bifurcation points on the primary solution path. Note that, the numerical solutions $u(x)$ on the primary path are all reflection symmetry-preserving. In fact, the numerical bifurcation paths even with coarser FE discretizations agree well with (\ref{eq.u_0_star}) except near the limit point. Through this work, we bring to attention that while numerical approaches are routinely used for parametrized nonlinear problems, a systematic verification of the numerical bifurcations cannot be discounted.\\  
\\
{\bf{Acknowledgments}}: The work of SSP was supported by IIT Bombay's IRCC seed grant RD/0524-IRCCSH0-010.

\bibliographystyle{plain}
\bibliography{lib}

@book{govaerts2000numerical,
  title={Numerical methods for bifurcations of dynamical equilibria},
  author={Govaerts, W. JF},
  year={2000},
  publisher={SIAM}
}

@article{kuiper1988spurious,
  title={On spurious numerical solutions for nonlinear eigenvalue problems},
  author={Kuiper, H. J},
  journal={The Rocky Mountain journal of mathematics},
  volume={18},
  number={2},
  pages={357--385},
  year={1988},
  publisher={JSTOR}
}

@article{liouville1853equation,
  title={Sur l'{\'e}quation aux diff{\'e}rences partielles},
  author={Liouville, J.},
  journal={Journal de math{\'e}matiques pures et appliqu{\'e}es},
  volume={18},
  pages={71--72},
  year={1853}
}

@book{peitgen1980nonlinear,
  title={Nonlinear elliptic boundary value problems versus their finite difference approximations},
  author={Peitgen, H.-O. and Saupe, D. and Schmitt, K.},
  year={1980}
}

@article{domokos2003nonlinear,
  title={On nonlinear boundary-value problems: ghosts, parasites and discretizations},
  author={Domokos, G. and Holmes, P.},
  journal={Proceedings of the Royal Society of London. Series A: Mathematical, Physical and Engineering Sciences},
  volume={459},
  number={2034},
  pages={1535--1561},
  year={2003},
  publisher={The Royal Society}
}

@incollection{ALLGOWER19751,
title = {On a Discretization of $y'' + \lambda y^{k} = 0$},
editor = {Miller, J. J. H.},
booktitle = {Topics in Numerical Analysis II},
publisher = {Academic Press},
pages = {1-15},
year = {1975},
isbn = {978-0-12-496952-0},
doi = {https://doi.org/10.1016/B978-0-12-496952-0.50008-8},
url = {https://www.sciencedirect.com/science/article/pii/B9780124969520500088},
author = {Allgower, E. L.}
}

@article{beyn1982spurious,
  title={Spurious Solutions for Discrete Superlinear Boundary Value Problems.},
  author={Beyn, W.-J. and Lorenz, J.},
  journal={Computing},
  volume={28},
  number={1},
  pages={43--51},
  year={1982}
}

@article{korman2003accurate,
  title={An accurate computation of the global solution curve for the Gelfand problem through a two point approximation},
  author={Korman, P.},
  journal={Applied mathematics and computation},
  volume={139},
  number={2-3},
  pages={363--369},
  year={2003},
  publisher={Elsevier}
}

@book{keller1987lectures,
  title={Lectures on Numerical Methods in Bifurcation Problems},
  author={Keller, H. B.},
  series={Lectures on mathematics and physics},
  year={1987},
  publisher={Tata Institute of Fundamental Research}
}

@book{kielhofer2011,
  title={Bifurcation Theory: An Introduction with Applications to Partial Differential Equations},
  author={Kielh{\"o}fer, H.},
  series={Applied Mathematical Sciences},
  year={2011},
  publisher={Springer}
}

@article{JACOBSEN2002283,
title = {The Liouville–Bratu–Gelfand Problem for Radial Operators},
journal = {Journal of Differential Equations},
volume = {184},
number = {1},
pages = {283-298},
year = {2002},
issn = {0022-0396},
doi = {https://doi.org/10.1006/jdeq.2001.4151},
url = {https://www.sciencedirect.com/science/article/pii/S0022039601941518},
author = {Jacobsen, J. and Schmitt, K.}}

@BOOK{Frank-Kamenetskii,
  title={``Diffusion and Heat Exchange in Chemical Kinetics"},
  author={Frank-Kamenetskii, D. A.},
  year={1955},
  publisher={Princeton University Press},
  url={https://ieeexplore.ieee.org/document/9452603},}

@BOOK{Chandrashekhar,
  author={Chandrasekhar, S.},
  title={``An Introduction to the Study of Stellar Structure,"},
  year={1957},
  volume={},
  number={},
  pages={},
  doi={},
  publisher={Dover Publications
 Inc., New York, NY},}

@BOOK{Gavalas1968,
  author={Gavalas, G. R.},
  title={‘‘Nonlinear Differential Equations of Chemically Reacting Systems,’’},
  year={1968},
  publisher={Springer-Verlag, New York Inc., New York,},}

@article{Farrell2015,
url = {https://doi.org/10.1137/140984798},
title = {Deflation techniques for finding distinct solutions of nonlinear partial differential equations},
author = {Farrell, P. E. and Birkisson, {\'A}. and Funke, S. W.},
pages = {A2026--A2045},
volume = {37},
number = {4},
journal = {SIAM Journal on Scientific Computing},
doi = {doi.org/10.1137/140984798},
year = {2015},
}

@article{stakgold1971,
url = {https://doi.org/10.1137/1013063},
title = {Branching of solutions of nonlinear equations},
author = {Stakgold, I.},
pages = {289-332},
volume = {13},
number = {3},
journal = {SIAM Review},
doi = {doi.org/10.1137/1013063},
year = {1971},
}

@article{rabinowitz1971,
url = {https://doi.org/10.1016/0022-1236(71)90030-9},
title = {Some global results for nonlinear eigenvalue problems},
author = {Rabinowitz, P.},
pages = {487-513},
volume = {7},
number = {3},
journal = {Journal of Functional Analysis},
doi = {doi.org/10.1016/0022-1236(71)90030-9},
year = {1971},
}

\end{document}